\documentclass[12pt,draftcls,onecolumn]{IEEEtran}

\usepackage{amssymb}
\usepackage{makeidx}         
\usepackage{amsmath}
\usepackage{multicol}        
\usepackage[bottom]{footmisc}
\usepackage{graphicx}
\usepackage{color}
\usepackage{amsfonts}
\usepackage[latin1]{inputenc}
\usepackage{epstopdf}
\usepackage[usenames,dvipsnames]{xcolor}
\usepackage[ruled,vlined]{algorithm2e}

\newcommand{\eq}{\begin{equation}}
\newcommand{\eeq}{\end{equation}}
\newcommand{\eqn}{\begin{eqnarray}}
\newcommand{\eeqn}{\end{eqnarray}}

\newcommand{\bsea}{\begin{subeqnarray}}
\newcommand{\esea}{\end{subeqnarray}}
\newcommand{\nn}{\nonumber}

\newcommand{\alg}[1]{\begin{align}  #1 \end{align}}

\def\bmat{\left[ \begin{array}}
\def\emat{\end{array} \right]}

\newcommand{\tr}{\mathop{\rm tr}}  

\newcommand{\Bc}{ \mathcal{B}}

\newcommand{\Dc}{ \mathcal{D}}

\newcommand{\Gc}{ \mathcal{G}}

\newcommand{\Nc}{ \mathcal{N}}

\newcommand{\Qc}{ \mathcal{Q}}

\newcommand{\Ds}{ \mathbb{D}}
\newcommand{\Es}{ \mathbb{E}}

\newcommand{\Rs}{ \mathbb{R}}

\newcounter{pippo}

\newtheorem{remark}{Remark}[section]
\newtheorem{teor}{Theorem}[section]
\newtheorem{corr}{Corollary}[section]
\newtheorem{propo}{Proposition}[section]
\newtheorem{lemm}{Lemma}[section]
\newtheorem{exam}{Example}
\newtheorem{probl}[pippo]{Problem}
\newtheorem{defn}{Definition}[section]

\newcommand{\teo}{\begin{teor}}
\newcommand{\eteo}{\end{teor}}
\newcommand{\cor}{\begin{corr}}
\newcommand{\ecor}{\end{corr}}
\newcommand{\prop}{\begin{propo}}
\newcommand{\eprop}{\end{propo}}
\newcommand{\lem}{\begin{lemm}}
\newcommand{\elem}{\end{lemm}}
\newcommand{\ex}{\begin{exam}}
\newcommand{\eex}{\end{exam}}
\newcommand{\pb}{\begin{probl}}
\newcommand{\epb}{\end{probl}}
\newcommand{\df}{\begin{defn}}
\newcommand{\edf}{\end{defn}}
\newcommand{\aprop}{\begin{apropo}}
\newcommand{\eaprop}{\end{apropo}}
\newcommand{\alem}{\begin{alemm}}
\newcommand{\ealem}{\end{alemm}}
\newcommand{\rem}{\begin{remark}}
\newcommand{\erem}{\end{remark}}

\newcommand{\Dd}  {\mathrm{D}}

\begin{document}

\title{Robust Kalman Filtering under Model Perturbations \thanks{This work has been partially supported by the FIRB project ``Learning meets time'' (RBFR12M3AC) funded by MIUR.}}

\author{Mattia~Zorzi \thanks{M. Zorzi is with the
Dipartimento di Ingegneria dell'Informazione, Universit\`a degli studi di
Padova, via Gradenigo 6/B, 35131 Padova, Italy
({\tt\small zorzimat@dei.unipd.it}).}}

\markboth{DRAFT}{Shell \MakeLowercase{\textit{et al.}}: Bare Demo of IEEEtran.cls for Journals}

\maketitle

\begin{abstract}
We consider a family of divergence-based minimax approaches to perform robust filtering. The mismodeling budget, or tolerance, is specified at each time increment of the model. 
More precisely, all possible model increments belong to a ball which is formed by placing a bound on the Tau-divergence family between the actual and the nominal model increment. 
Then, the robust filter is obtained by minimizing the mean square error according to the least favorable model in that ball. It turns out that the solution is a family of Kalman like filters. Their gain matrix is updated according to a risk sensitive like iteration where the risk sensitivity parameter is now time varying. As a consequence, we also extend the risk sensitive filter to a family of risk sensitive like filters according to the Tau-divergence family. 

\end{abstract}

\begin{IEEEkeywords}
Robust Kalman filtering, Tau-divergence family, minimax problem, risk sensitive filtering.
\end{IEEEkeywords}

\section{Introduction}

Kalman filter is ubiquitous in many applications. The main reason is due by its iterative structure, allowing its implementation very simple. On the other hand, this filter is designed with respect to a linear state space model. The latter is often inadequate to describe phenomena, accordingly the resulting Kalman filter does not perform well in the practice. Since the beginning, therefore, it was clear the importance to develop robust versions of the standard Kalman filter.

Robust filtering can be performed according to the risk sensitive approach, \cite{RISK_WHITTLE_1980,OPTIMAL_SPEYER_FAN_BANAVAR_1992,RISK_PROP_BANAVAR_SPEIER_1998,H_INF_HASSIBI_SAYED_KAILATH_1999,speyer_2008,LEVY_ZORZI_RISK_CONTRACTION}. Here, the robust estimator is designed according to the nominal model 
but  in such a way to avoid large errors. The sensitivity to large errors is tuned by the so called risk sensitivity parameter. 
It is worth noting this approach has been interpreted as a minimax problem
\cite{boel2002robustness,YOON_2004,HANSEN_SARGENT_2005,HANSEN_SARGENT_2007,ROBUSTNESS_HANSENSARGENT_2008}. The appealing aspect of the risk sensitive approach is that the solution is a Kaman like filter.
On the other hand, the risk sensitive parameter is not explicitly connected to the discrepancy between the actual and the nominal model. 
Recently, a divergence-based minimax approach has been proposed in \cite{LEVY_NIKOUKHAH_2004,ROBUST_STATE_SPACE_LEVY_NIKOUKHAH_2013,ZORZI_CONTRACTION_CDC}. More precisely, in \cite{LEVY_NIKOUKHAH_2004}
the robust static estimation problem of a signal given noisy observations has been considered. Here, all possible models belong to a ball which is formed by placing a bound on the {\em Kullback-Leibler} divergence between the actual and the nominal model. 
This bound, say tolerance, represents the mismodeling budged. 
Then, the robust filter is obtained by minimizing the mean square error according to the least favorable model in this ball. 
It turns out that the Bayes estimator is robust under model uncertainty characterized by this ball. In \cite{ROBUST_STATE_SPACE_LEVY_NIKOUKHAH_2013}, a dynamic extension to this problem (i.e. a robust filtering problem) has been considered. More precisely, 
drawing inspiration from \cite{GHAOUI_CALAFIORE_2001,SAYED_FRAMEWORK_2001}, 
the mismodeling budged is specified to each time increment of the model, that is the model uncertainty is expressed in an incremental way. Roughly speaking, the idea is to iterate the Bayes estimator with the least favorable statistics found in \cite{LEVY_NIKOUKHAH_2004}. It turns out the robust estimator has a Kalman like structure. More precisely, it is a risk sensitivity like filter, where the risk sensitivity parameter is now time varying.

In \cite{OPTIMALITY_ZORZI}, the robust static estimation problem proposed in \cite{LEVY_NIKOUKHAH_2004} has been extended, in the Gaussian case, to a family of uncertainty classes. The latter is formed by placing a bound on a set of divergences (called $\tau$-divergence family,
\cite{OPTIMAL_PREDICTION_ZORZI_2014}) between the actual and the nominal model. This particular divergence family is chosen because, in contrast to the alpha and the beta family \cite{ALPHA,BETA,DUAL}, it allows to characterize uncertainty balls for which the Bayes estimator is still robust.

The contribution of this paper is
to extend the robust Kalman filter in \cite{ROBUST_STATE_SPACE_LEVY_NIKOUKHAH_2013}  to a family of robust Kalman filters 
parametrized by the $\tau$-divergence family using the results in \cite{OPTIMALITY_ZORZI}. This family of filters is characterized by a time varying risk sensitive parameter. Therefore, by adopting the perspective given in \cite{boel2002robustness}, we also extend the risk sensitive filter 
to a family of risk sensitive like  filters parametrized by the $\tau$-divergence family, say $\tau$-risk sensitive filters.
Finally, we present a simulation study which shows that parameter $\tau$ tunes how conservative the robust filter is.

In the paper we will use the following notation. $\|x\|$ denotes the Euclidean norm of $x\in\Rs^n$. $\|x\|_A$ denotes the weighted Euclidean norm with $A$ symmetric and positive definite. 

\section{Robust Static Estimation} \label{sec_static}
We review the robust static estimation problem under model uncertainty characterized by the $\tau$-divergence family introduced in \cite{OPTIMALITY_ZORZI}.
Let $x\in\Rs^n$ and $y\in\Rs^p$ be two jointly Gaussian random vectors.
Let  $z:= \left[
       \;             x^T \;y^T \;
         \right]^T$.  Its joint nominal probability  density $f$ is
\alg{ \small f(z)= \frac{1}{\sqrt{(2\pi)^{p+n} \det K_z }}\exp\left(-\frac{1}{2}(z-m_z)^T K_z^{-1} (z-m_z)\right) \nn}
where the mean vector $m_z\in\Rs^{n+p}$ and the covariance matrix $K_z\in\Qc_+^{n+p}$ are known. We conformably partition the mean vector and the covariance matrix of $z$ according to $x$
and $y$:
\alg{ m_z=\left[
          \begin{array}{c}
            m_x \\
            m_y \\
          \end{array}
        \right], \;\; K_z=\left[
                            \begin{array}{cc}
                              K_x & K_{xy} \\
                              K_{yx} & K_y\\
                            \end{array}
                          \right]\nn
.} Let $\tilde f$ denote the actual joint probability density of $z$
\alg{ \tilde f(z)= \frac{1}{\sqrt{(2\pi)^{p+n} \det \tilde K_z }}\exp\left(-\frac{1}{2}(z-m_z)^T \tilde K_z^{-1} (z-m_z)\right) \nn}
where the mean vector $\tilde m_z\in\Rs^{n+p}$ and the covariance matrix $\tilde K_z\in\Qc_+^{n+p}$ are unknown.
Since both $f$ and $\tilde f$ are Gaussian, the deviation between $f$ and $ \tilde f$ may be directly measured by the
deviation between $(m_z,K_z)$ and $(\tilde m_z,\tilde K_z$) through the $\tau$-divergence, \cite{OPTIMALITY_ZORZI}:
\alg{& \Dc_\tau(\tilde f\| f)=\nn\\
& \left\{
                               \begin{array}{ll}
                                 \|\Delta m_z\|^2_{K_Z^{-1}}+\tr\left(-\log(\tilde K_z K_z^{-1})\right. &\\ \hspace{0.3cm}\left.+\tilde K_z K_z^{-1}-I_{n+p} \right), & \hbox{$\tau=0$} \\
                                 \|\Delta m_z\|^2_{\frac{1}{1-\tau}K_Z^{-1}}+\tr\left(-\frac{1}{\tau(1-\tau)}(L_z^{-1}\tilde K_z L_z^{-T})^{\tau}\right. & \\ \hspace{0.3cm}\left.  +\frac{1}{1-\tau} \tilde K_z K_z^{-1}+\frac{1}{\tau}I_{n+p}\right), & \hspace{-0.4cm}\hbox{$0<\tau<1$} \\
                               \delta_\infty(\Delta m_z)+  \tr\left( L_z^{-1}\tilde K_z L_z^{-T}\log(L_z^{-1}\tilde K_z                                    L_z^{-T})\right. &\\ \label{def_tau_div}
                               \hspace{0.3cm}\left.-\tilde K_z K_z^{-1}+I_{n+p}\right), & \hbox{$\tau=1$}
                               \end{array}
                             \right. }
where $L_z$ is such that $K_z=L_zL_z^T$, $\Delta m_z=\tilde m_z - m_z$, and
\alg{ \delta_\infty(v):=\left\{\begin{array}{cc} 0 &  \hbox{if }v=0  \\ \infty &  \hbox{otherwise}.\\ \end{array}\right. \nn}
 Note that, $\Dc_\tau(\tilde f \| f)\geq 0$ and equality holds if and only if $\tilde f=f$. This divergence takes root in the prediction theory. Let $e^N=L_z^{-1}(z-m_z)$ with $z\sim \tilde f$. $e^N$ can be understood  as a normalized prediction error, where $m_z$ represents the minimum variance prediction of $z$ based on $f$. If $\tilde f=f$, then $e^N$ has zero mean and covariance matrix $I$.
Hence, this divergence measures the discrepancy between $e^N$ and the Gaussian random vector with zero mean and covariance $I$. We consider the closed ball centered on $f$:
\eq \label{ball_static}\Bc_\tau:=\{ f \hbox{ s.t. } \Dc_\tau(\tilde f\| f)\leq c\}\eeq
where $c\in\Rs_+$ is a fixed tolerance. Accordingly, $\Bc_\tau$ represents the set of all possible probability densities of $z$ consistent with the allowed mismodelling budget.

The robust estimator of $x$ given $y$ is designed according to the
minimax point of view \cite{LEVY_NIKOUKHAH_2004,ROBUSTNESS_HANSENSARGENT_2008}. More precisely, whenever
we seek to design an estimator minimizing a suitable loss function, an hostile player, say ``nature'', conspires to select the worst possible probability density in $\Bc_\tau$. 
Let $g(y)$ denote an estimator of $x$ based on the observation vector $y$.
The optimal robust filter is solution to the following minimax problem
 \eq \label{minimax_problem_fg}\underset{g\in\Gc}{\min}\;\underset{\tilde f\in\Bc_\tau}{\max} \;J(\tilde f,g)\eeq
where
\alg{ \label{Formula_J} J(\tilde f,g)&= \Es_{\tilde f}[\| x-g(y)\|^2]=\int _{\Rs^{n+p}} \| x-g(y)\|^2 \tilde f(z)\mathrm{d} z\nn }
denotes the mean square error and $\Gc$ denotes the set of all estimators $g(y)$ such that $\Es_{\tilde f}[\|g(y)\|^2]$ is finite for any $\tilde f \in\Bc_\tau$.

\teo \label{theorem_static}Let $0\leq \tau\leq 1$.  The least favorable probability density $\tilde f^\circ$ has mean vector $\tilde m_z^\circ=m_z$ and covariance matrix
with the following structure \eq \label{def_K_tilda}\tilde K^\circ_z=\left[
                                                                         \begin{array}{cc}
                                                                           \tilde K_x & K_{xy} \\
                                                                           K_{yx} &  K_y \\
                                                                         \end{array}
                                                                       \right]
\eeq wherein only the covariance of $x$ is perturbed with respect to the nominal covariance matrix.
Let \alg{ \label{def_delle_P}  P &=  K_x-K_{xy}K_y^{-1}K_{yx}\nn\\
V&= \tilde K_x-K_{xy}K_y^{-1}K_{yx}}
denote the nominal and the perturbed {\em a posteriori} covariance matrix of $x$ given $y$. Then,
\eq \label{def_V} V=\left\{
                                   \begin{array}{ll}
                                     L_P\left(I_{n}-\theta(1-\tau)L_P^TL_P\right)^{\frac{1}{\tau-1}}L_P^T, & 0 \leq \tau<1 \\
                  L_P\exp\left(\theta L_P^TL_P\right) L_P^T                   , & \tau=1
                                   \end{array}
                                 \right.\eeq
where $L_P$ is such that $P=L_PL_P^T$. Here $\theta^{-1}$, with $\theta^{-1}>(1-\tau)\| P\|$, is the unique {\em Lagrange} multiplier
satisfying the relation $c=\gamma_\tau(P,\theta)$
where
\alg{ \label{def_gamma} & \gamma_\tau(P,\theta)=\nn\\ & \hspace{0.1cm}\left\{
      \begin{array}{ll}
       -\log\det(I_n-\theta P)^{-1}+\tr((I_n-\theta P)^{-1}-I_n) , &  \tau=0\\
\tr (-\frac{1}{\tau(1-\tau)}(I_n-\theta(1-\tau) L_{P}^TL_{P})^{\frac{\tau}{\tau-1}}  \\
 \hspace{0.3cm}+\frac{1}{1-\tau} (I_n-\theta(1-\tau) L_{P}^TL_{P})^{\frac{1}{\tau-1}}+\frac{1}{\tau}I_n)        , & \hspace{-0.6cm}0 <\tau<1\\
      \tr(\exp(\theta L_{P}^T L_{P})(\theta L_{P}^T L_{P}-I_n)+I_n)  , & \tau=1
      \end{array}
    \right.}

 The optimal robust estimator is the Bayes estimator \eq \label{Opt_estimator}g^\circ(y)=G^\circ(y-m_y)+m_x\eeq
with $G^\circ=K_{xy}K_y^{-1}$.
\eteo

Theorem \ref{theorem_static} shows that the Bayes estimator is robust with respect to the uncertainty class, parametrized by $\tau$, in 
(\ref{ball_static}). Clearly, this optimality holds in the Gaussian case. Without this assumption, the least favorable probability density could be more different than the one in Theorem \ref{theorem_static}.

\cor \label{corollario_static} Let $\theta>0$ be {\em  a priori} fixed and such that $\theta^{-1}>(1-\tau)\|P\|$.
Consider the minimax problem
 \alg{ \underset{g\in\Gc}{\min}\;\underset{\tilde f\in\Bc_\tau}{\max} \; \Es_{\tilde f}[\|x-g(y)\|^2]-\theta^{-1}\Dc_\tau(\tilde f\| f)\nn}
 where $\Bc_\tau=\{\tilde f \hbox{ s.t. } \Dc_\tau(\tilde f \| f)< \infty \}$ and $\Gc$ is the set of all estimators such that $\Es_{\tilde f}[\| g(y)\|^2]$
is finite for any $\tilde f \in\Bc_\tau$. Then, the least favorable probability density $\tilde f_t ^\circ$ has mean vector $\tilde m_z^\circ =m_z$ and covariance matrix $\tilde K_z^\circ$
as in (\ref{def_K_tilda}). The perturbed {\em a posteriori} covariance matrix $V$ of $x$ given $y$ is  in (\ref{def_delle_P}). Moreover, its relation with $P$ is given by (\ref{def_V}) where $\theta$ now has been {\em a priori} chosen. The optimal estimator is the Bayes estimator (\ref{Opt_estimator}).\ecor

\section{Robust filtering problem} \label{sec_pb_formulation}
We consider a nominal Gauss-Markov state space model of the form
\alg{ \label{state_space_model}x_{t+1}=A_tx_t+B_t v_t\nn\\
y_t=C_tx_t+D_t v_t}
where $x_t\in\Rs^n$ is the state process, $y_t\in\Rs^p$ is the observation process, and $v_t\in\Rs^m$ is WGN with unit variance, i.e.
$ \Es[v_t v_s^T]=I_m\delta_{t-s} $
where $\delta_{t}$ denotes the Kronecker delta function.
We assume that the noise $v_t$ is independent of the initial state, whose nominal distribution is given by
$ f_0(x_0)\sim \Nc(\hat x_0,V_0)$.
Let $ z_t=\left[\;
                x_{t+1}^T \;
                y_{t}^T \;
            \right]^T
.$
Model (\ref{state_space_model}) is characterized by the nominal transition probability density of $z_t$
given $x_t$:
\eq \label{nominal_phi}\phi_t(z_t|x_t)\sim\Nc\left(\left[
                             \begin{array}{c}
                               A_t \\
                               C_t \\
                             \end{array}
                           \right]x_t, \left[
                             \begin{array}{c}
                               B_t \\
                               D_t \\
                             \end{array}
                           \right] \left[
                                     \begin{array}{cc}
                                       B_t^T & D_t^T \\
                                     \end{array}
                                   \right]
\right).\eeq
As noticed in \cite{ROBUST_STATE_SPACE_LEVY_NIKOUKHAH_2013}, when entropy-like indexes are used to measure the proximity of statistical models, all the relations between dynamic variables or observations
should be uncertain,
otherwise those indexes take infinite value. To avoid such  a situation, we assume that the noise $v_t$ affects all the components of the dynamics and observations in (\ref{state_space_model}), possibly with a very small variance for relations which are viewed as almost certain. Therefore, the covariance matrix
\alg{ K_{z_t|x_t}=\left[
      \begin{array}{c}
        B_t \\
        D_t \\
      \end{array}
    \right]\left[
             \begin{array}{cc}
               B_t^T & D_t^T \\
             \end{array}
           \right]
\nn}
is positive definite. Moreover, the matrix
$ \Gamma_t=\left[\;
        B_t^T\;
        D_t^T \;
    \right]^T$
has full column rank, and without loss of generality we can assume $\Gamma_t$ is square and invertible, so that $m=n+p$. Otherwise, we can compress the column space of $\Gamma_t$ and remove noise components which do not affect model (\ref{state_space_model}).

We adopt the minimax approach proposed \cite[Section III]{ROBUST_STATE_SPACE_LEVY_NIKOUKHAH_2013} to characterize the robust filter. Let $\tilde \phi_s(z_t|x_t)$ be the least favorable transition probability density of $z_t$ given $x_t$. 
Let $\tilde f_t(x_t|Y_{t-1})$ be the {\em a priori} probability density of $x_t$ conditioned on the observations $Y_{t-1}=\{ y_s,\;\; 0\leq s\leq t-1\}$ and based on the least favorable model. We introduce the marginal probability densities
\alg{ \label{def_nominal_marg_density}\bar f_t(z_t|Y_{t-1})&=\int \phi_t(z_t|x_t)\tilde f_t(x_t|Y_{t-1})\mathrm{d}x_t\\
\label{def_leastfav_marg_density}\tilde f_t(z_t|Y_{t-1})&=\int \tilde \phi_t(z_t|x_t)\tilde f_t(x_t|Y_{t-1})\mathrm{d}x_t.}
 Note that, $\bar f_t(z_t|Y_{t-1})$ can be viewed as the
  pseudo-nominal density of $z_t$ conditioned on $Y_{t-1}$ computed from the conditional least favorable  density $\tilde f_t(x_t|Y_{t-1})$ and the nominal
transition probability density $\phi_t(x_t|z_t)$.
As in \cite{ROBUST_STATE_SPACE_LEVY_NIKOUKHAH_2013}, we assume that
\eq  \label{def_a_priori_least_fav}\tilde f_t(x_t| Y_{t-1})\sim \Nc(\hat x_t,V_t).\eeq
In this way the conditional probability density $\bar f_t(z_t|Y_{t-1})$ is Gaussian. We make the additional assumption that $\tilde \phi(z_t|x_t)$ is such that $\tilde f_t(z_t|Y_{t-1})$ is Gaussian. In \cite{ROBUST_STATE_SPACE_LEVY_NIKOUKHAH_2013}, the latter assumption was not made. However, it is worth noting that the least favorable solution found is such that (\ref{def_leastfav_marg_density}) is Gaussian, see Remark \ref{remark_lavoro_prec}.
Therefore, we can measure the deviance between $\tilde \phi_t$ and $\phi_t$ as deviance between $\tilde f_t(z_t|Y_{t-1})$
and $\bar f_t(z_t|Y_{t-1})$ using the $\tau$-divergence (\ref{def_tau_div}). Then, we assume that $\tilde \phi_t$
belongs to the closed ball about $\phi_t$:
$ \Bc_{t,\tau}=\{\tilde \phi_t(z_t|x_t)\hbox{ s.t. } \Dc_\tau(\tilde f_t\| \bar f_t)\leq c_t\}$
where $c_t\in\Rs_+$ is the tolerance specified at each time step.
Let $\Gc_t$ denote the class of estimators with finite second-order moments with respect all densities $\tilde \phi_t(z_t|x_t) \tilde f_t(x_t|Y_{t-1})$ such that $\tilde \phi_t(z_t|x_t)\in\Bc_{t,\tau}$. The the robust filter is characterized by the following minimax problem
\eq \label{minimax_game} (g_t^\circ,\tilde \phi_t ^\circ)=\arg \underset{g_t\in\Gc_t}{\min}\underset{\tilde \phi_t\in\Bc_{t,\tau}}{\max} J_t(\tilde \phi_t,g_t)\eeq
where \alg{ J_t & (\tilde \phi_t,g_t)= \Es_{\tilde f_t}[\|x_{t+1}-g_t(y_t)\|^2|Y_{t-1}]\nn\\
 &=\int\int \|x_{t+1}-g_t(y_t)\|^2 \tilde \phi_t(z_t|x_t) \tilde f_t(x_t| Y_{t-1}) \mathrm{d} x_t \mathrm{d}z_t \nn } denotes the mean square error of the estimator $\hat x_{t+1}=g_t^\circ(y_t)$
of $x_{t+1}$ evaluated with respect to the transition density $\tilde \phi_t$ in $\Bc_{t,\tau}$. It is
worth noting that $\hat x_{t+1}$ depends on $Y_t$, and not only on $y_t$, but this dependency is suppressed to simplify the notations.

\rem In the minimax problem (\ref{minimax_game}) we require that $\tilde f_t(z_t|Y_{t-1})$ defined in (\ref{def_leastfav_marg_density}) is a conditional probability
density, that is
\alg{\label{cond_normalizzazione}  \int\int \tilde \phi_t(z_t|x_t)\tilde f_t(x_t|Y_{t-1})\mathrm{d}z_t\mathrm{d}x_t=1,} but we do not require that $\tilde \phi_t(z_t|x_t)$ is a 
transition probability density for each $x_t$.
Therefore, the {\em a priori} conditional probability density $\tilde f_t(x_t|Y_{t-1})$ is not required to coincide with the
{\em a posteriori} one computed from $\tilde \phi_t(z_t|x_t)\tilde f_t(x_t|Y_{t-1})$. \erem

\rem \label{remark_lavoro_prec} In \cite{ROBUST_STATE_SPACE_LEVY_NIKOUKHAH_2013}, $\tilde \phi_t(z_t|x_t)$ is not required to be such that (\ref{def_leastfav_marg_density}) is Gaussian. The constraint on $\tilde \phi_t(z_t|x_t)$ is that $\Ds_{KL}(\tilde f_t \| \bar f_t)\leq c_t$ where $\Ds_{KL}$ is the {\em Kullback-Leibler} divergence among probability densities.
 On the other hand, the solution $\tilde \phi_t^\circ(z_t|x_t)$ to the corresponding 
minimax problem is such that (\ref{def_leastfav_marg_density}) is Gaussian, see 
\cite[Formula (16)]{ROBUST_STATE_SPACE_LEVY_NIKOUKHAH_2013}. Hence, the corresponding $\tilde f_t^\circ(z_t|Y_{t-1})$ is Gaussian. Note that, $\Dd_{KL}(\tilde f_t\| \bar f_t)=\Dc_{0}(\tilde f_t\| \bar f_t)$ when $\tilde f_t$, $\bar f_t$ are Gaussian. We conclude that, for  $\tau=0$, the solution to (\ref{minimax_game}) coincides with the one in \cite{ROBUST_STATE_SPACE_LEVY_NIKOUKHAH_2013}. 
 \erem

\section{Robust Kalman Filters} \label{sec_rob_kalman}
We show that the optimal robust estimator solution to the minimax problem (\ref{minimax_game}) is a Kalman like filter 
parametrized by $\tau$. In this way, we obtain a parametric family of robust {Kalman} filters. 
First, Problem (\ref{minimax_game}) can be reformulated as the static minimax problem
(\ref{minimax_problem_fg}). Consider the ball
$ \bar\Bc_{t,\tau}=\{\tilde f_t(z_t|Y_{t-1}) \hbox{ s.t. } \Dc_\tau(\tilde f_t\| \bar f_t )\leq c_t\}$
which is the set of all probability densities having structure (\ref{def_leastfav_marg_density}) with $\tilde \phi_t\in\Bc_{t,\tau}$. The equivalent minimax problem is
\alg{ (\tilde f_t^\circ ,g^\circ_t)= \arg \underset{g_t\in\Gc_t}{\min}\underset{\tilde f_t\in \bar \Bc_{t,\tau}}{\max} \bar J_t(\tilde f_t,g_t)\nn}
where \alg{ \bar J_t(\tilde f_t,g_t)=\int \|x_{t+1}-g_t(y_t)\|^2 \tilde f_t(z_t|Y_{t-1})\mathrm{d}z_t.\nn}
In view of (\ref{nominal_phi}) and (\ref{def_a_priori_least_fav}), the pseudo-nominal density is Gaussian
\eq \label{bar_f_t_Gaussian}\bar f_t(z_t| Y_{t-1})\sim \Nc\left( \left[
                                      \begin{array}{c}
                                        A_t \\
                                        C_t \\
                                      \end{array}
                                    \right]\hat x_t, K_{z_t}
\right)\eeq where the conditional covariance matrix $K_{z_t}$ is given by
\alg{ \label{def_K_zt}K_{z_t}&=\left[\begin{array}{cc} K_{x_{t+1}} & K_{x_{t+1},y_t}  \\ K_{y_t,x_{t+1}} &  K_{y_t}\end{array}\right]\nn\\ &= \left[
                                      \begin{array}{c}
                                        A_t \\
                                        C_t \\
                                      \end{array}
                                    \right] V_t \left[
                                                  \begin{array}{cc}
                                                    A_t^T & C_t^T  \\
                                                  \end{array}
                                                \right]+\left[
                                      \begin{array}{c}
                                        B_t \\
                                        D_t \\
                                      \end{array}
                                    \right] \left[
                                                  \begin{array}{cc}
                                                    B_t^T & D_t^T  \\
                                                  \end{array}
                                                \right].
                                    }

Applying Theorem \ref{theorem_static}
with $f \mapsto \bar f_t $, $\tilde f  \mapsto \tilde f_t $
and $g \mapsto g_t $, the least favorable conditional density $\tilde f_t^\circ(z_t|Y_{t-1})$ is such that
\eq \label{def_tilda_f_t} \tilde f_t^\circ(z_t|Y_{t-1})\sim \Nc\left( \left[
                                                  \begin{array}{c}
                                                    A_t \\
                                                    C_t \\
                                                  \end{array}
                                                \right] \hat x_t, \tilde K_{z_t}^\circ
\right)\eeq where the least favorable conditional covariance matrix is \alg{ \tilde K_{z_t}^\circ=\left[
                           \begin{array}{cc}
                            \tilde  K_{x_{t+1}} & K_{ x_{x+1}y_t} \\
                             K_{y_{t} x_{t+1}} & K_{y_t} \\
                           \end{array}
                         \right].\nn
} 
Let
\alg{ P_{t+1}&=K_{x_{t+1}}-K_{x_{t+1},y_t}K_{y_t}^{-1}K_{y_t, x_{t+1}}\nn \\
 V_{t+1}&= \tilde K_{x_{t+1}}-K_{x_{t+1},y_t}K_{y_t}^{-1}K_{y_t, x_{t+1}}\nn }
denote the nominal and the least favorable conditional covariance of $x_{x+1}$ given $Y_t$. Then,
\alg{& V_{t+1}=\nn\\ &\left\{
      \begin{array}{ll}
        L_{P_{t+1}}\left(I_n-\theta_t(1-\tau) L_{P_{t+1}}^TL_{P_{t+1}}\right)^{\frac{1}{\tau-1}} L_{P_{t+1}}^T , & \hspace{-0.1cm}  0<\tau<1\\
L_{P_{t+1}}\exp\left(\theta_t L_{P_{t+1}}^T L_{P_{t+1}}\right) L_{P_{t+1}}^T, & \tau=1
      \end{array}
    \right. \nn 
}
 where $L_{P_{t+1}}$ is such that $P_{t+1}=L_{P_{t+1}}L_{P_{t+1}}^T$ and $\theta_t^{-1}>(1-\tau)\|P_{t+1}\|$
is the unique solution to $c_t=\gamma_\tau(P_{t+1},\theta_t)$ where
$\gamma_\tau$ has been defined in (\ref{def_gamma}).
The optimal robust estimator takes the form,
\eq \label{rob_est_formula1} \hat x_{t+1}=g_t^\circ(y_t)=A_t \hat x_t+G_t(y_t-C_t \hat x_t) \eeq
with gain matrix $G_t=K_{x_{t+1},y_t} K_{y_t}^{-1}$.
From (\ref{def_K_zt}), we obtain
\alg{ \label{rob_est_formula2} G_t&=(A_t V_t C_t^T+B_t D_t^T)(C_t V_t C_t^T+D_t D_t^T)^{-1} \nn\\
P_{t+1}&= A_tV_t A_t^T- G_t(C_t V_t C_t^T+D_t D_t^T) G_t^T+B_t B_t^T.}

Algorithm \ref{algorithm} shows the iterative scheme of the optimal robust estimator  we found for the case $0\leq \tau<1$. 
The algorithm for the limit case $\tau=1$ is the same  with the exception that the updating of $V_{t+1}$ is different.
It is clear that the robust filter has the same iterative structure of
 the Kalman filter with the exception that $P_t$ is applied a distortion through matrix $V_t$. 
 In particular, $G_t$ is governed by a Riccati-like equation. 
\begin{algorithm}\label{algorithm}
 \LinesNumbered
 \SetKwInOut{Input}{Input}\SetKwInOut{Output}{Output}
 \Input{$c_t$, $\hat x_t$, $V_t$, $y_t$}
 \Output{$\hat x_{t+1}$, $V_{t+1}$}
 $G_t=(A_t V_t C_t+B_tD_t^T)^T(C_t V_t C_t^T+D_t D_t^T)^{-1}$\\
$\hat x_{t+1}=A_t \hat x_t+G_t(y_t-C_t \hat x_t)$\\
 $P_{t+1}= A_tV_t A_t^T- G_t(C_t V_t C_t^T+D_t D_t^T) G_t^T+B_t B_t^T $\\
Find $\theta_t$ such  that $c_t=\gamma_{\tau}(P_{t+1},\theta_t)$\\
Compute $V_{t+1}=L_{P_{t+1}}\left(I_n-\theta_t(1-\tau) L_{P_{t+1}}^TL_{P_{t+1}}\right)^{\frac{1}{\tau-1}} L_{P_{t+1}}^T$
 \caption{Robust Kalman filter at time $t$}
\end{algorithm}

It remains to characterize the least favorable transition density $\tilde \phi^\circ_t (z_t|x_t)$. It is not
difficult to prove that, \cite[Theorem 2.1]{OPTIMALITY_ZORZI},
\alg{ K_{z_t}^{-1}&=\left[
       \begin{array}{cc}
         I_{n} & 0 \\
         -G_t^T & I_{p} \\
       \end{array}
     \right] \left[
       \begin{array}{cc}
         P_{t+1}^{-1} & 0 \\
         0 & K_{y_t}^{-1} \\
       \end{array}
     \right] \left[
       \begin{array}{cc}
         I_{n} & -G_t \\
         0 & I_{p} \\
       \end{array}
     \right]\nn\\
(\tilde K_{z_t}^\circ)^{-1}&=\left[
       \begin{array}{cc}
         I_{n} & 0 \\
         -G_t^T & I_{p} \\
       \end{array}
     \right] \left[
       \begin{array}{cc}
         V_{t+1}^{-1} & 0 \\
         0 & K_{y_t}^{-1} \\
       \end{array}
     \right] \left[
       \begin{array}{cc}
         I_{n} & -G_t \\
         0 & I_{p} \\
       \end{array}
     \right] .\nn }
Accordingly
\alg{ K_{z_t}^{-1}-(\tilde K_{z_t}^\circ)^{-1}= \left[
       \begin{array}{c}
         I_{n}  \\
         -G_t^T \\
       \end{array}
     \right]\Phi_{t} \left[
       \begin{array}{cc}
         I_{n} & -G_t \\
       \end{array}
     \right]\nn }
where $\Phi_{t}=P_{t+1}^{-1}-V_{t+1}^{-1}$ which is positive definite.
Let $e_t=x_t-\hat x_t$ denote the estimation error. Define 
\alg{\hat m_{z_t}&=\Es_{f_t}[z_t|Y_{t-1}]=\Es_{\tilde f_t^\circ}[z_t|Y_{t-1}]=\left[
                                                                        \begin{array}{cc}
                                                                          A_t^T & C_t^T \\
                                                                        \end{array}
                                                                      \right]^T \hat x_t .\nn }
Therefore,
\alg{ \label{relation_z_et_Q_t+1}  (z_t & -\hat m_{z_t})^T  ( K_{z_t}^{-1}-(\tilde K_{z_t}^\circ)^{-1})(z_t-\hat m_{z_t})\nn\\
& =
(z_t-\hat m_{z_t})^T \left[
       \begin{array}{c}
         I_{n}  \\
         -G_t^T \\
       \end{array}
     \right]\Phi_{t} \left[
       \begin{array}{cc}
         I_{n} & -G_t \\
       \end{array}
     \right] (z_t-\hat m_{z_t})\nn\\
& =\| x_{t+1}-(A_t \hat x_t+G_t(y_t-C_t \hat x_t)) \|_{\Phi_{t}}^2
\nn\\ &=\| x_{t+1}-\hat x_{t+1} \|_{\Phi_{t}}^2 =\|e_{t+1}\|_{\Phi_{t}}^2 .\nn }
By (\ref{def_tilda_f_t}) and (\ref{bar_f_t_Gaussian}), we have
\alg{ \tilde f_t^\circ(z_t|Y_{t-1})
 \sim   \exp\left(\frac{1}{2} \|e_{t+1}\|^2_{\Phi_{t}}\right)\bar f_t(z_t|Y_{t-1}).\nn}
By (\ref{def_nominal_marg_density}), we obtain
\alg {\tilde f_t^\circ(z_t|Y_{t-1})  \sim  \int  \exp\left(\frac{1}{2} \|e_{t+1}\|^2_{\Phi_{t}}\right)\phi_t(z_t|x_t) \tilde f_t(x_t|Y_{t-1}) \mathrm{d}x_t. \nn  }
and by  (\ref{def_leastfav_marg_density}) we conclude that
\eq \label{def_tilde_phi_opt}\tilde \phi_t^\circ(z_t|x_t)=\frac{1}{M_t(\Phi_{t})} \exp\left(\frac{1}{2} \|e_{t+1}\|^2_{\Phi_{t}}\right)\phi_t(z_t|x_t) \eeq
where the normalizing constant $M_t(\Phi_{t})$ is such that (\ref{cond_normalizzazione}) holds.
It is worth noting that in the case $\tau=0$, i.e. the case considered in \cite{ROBUST_STATE_SPACE_LEVY_NIKOUKHAH_2013}, the distortion is a radial function of the estimation error $e_{t+1}$, because $\Phi_t=\theta_t^{-1} I_n$ for $\tau=0$.
On the contrary, in the case $\tau\neq 0$ such distortion is nonradial.

\section{Least-Favorable Model} \label{sec_least_fav}

For simulation and performance evaluation purposes, in particular for choosing parameters $c_t$ and $\tau$, it is important to characterize the least favorable model which is the solution to (\ref{minimax_game}).
The idea is to characterize it through (\ref{def_tilde_phi_opt}). Note that, there is  a one to one correspondence
between $z_t$ and $v_t$, given $x_t$, through the relation 
 \alg{ v_t=\Gamma_t^{-1} \left(z_t-\left[\begin{array}{c} A_t\\ C_t\end{array}\right]x_t\right)\nn }
because matrix $\Gamma_t$ is invertible. Accordingly, we can characterize the least favorable model
with model (\ref{state_space_model})
where the distortion has been moved now in noise  $v_t$. 
 Applying the same arguments used in \cite{ROBUST_STATE_SPACE_LEVY_NIKOUKHAH_2013}, see also \cite[Section 17.7]{ROBUSTNESS_HANSENSARGENT_2008},
it is not difficult to prove that the least favorable probability density of $v_t$ depends on $e_t$ and is distributed as 
$\tilde \psi_t(v_t|e_t)\sim \Nc(H_t e_t,\tilde K_{v_t})$
where
$\tilde K_{v_t}=(I_{n+p}-(B_t-G_tD_t)^T (\Omega_{t+1}^{-1}+\Phi_t) (B_t-G_tD_t))^{-1} $
and $H_t=\tilde K_{v_t} (B_t-G_tD_t)^T  (\Omega_{t+1}^{-1}+\Phi_t)(A_t-G_tC_t)$.
Matrix $\Omega_{t+1}^{-1}$ is computed from the backward recursion
\alg{\small \label{recursion_in_omega}\Omega_t^{-1}=(A_t & -G_tC_t)^T(\Omega_{t+1}^{-1}+\Phi_t) (A_t-G_tC_t) +H_t ^T \tilde K_{v_t}^{-1} H_t}
where the final point can be initialized with $\Omega_{T+1}^{-1}=0$ and $T$ is the simulation horizon.  The backward recursion is due by the fact that integrating $\tilde \phi_t^\circ(z_t|x_t)$
over $z_t$ we obtain a positive function of $e_t$, therefore the ``nature'' has the opportunity to change retroactively
the least favorable density of $x_t$. It is not difficult to see that the least favorable model admits a state space representation
with matrices \alg{\label{ls_model}\tilde A_t &= \left[\begin{array}{cc} A_t &  B_t H_t \\  0&  A_t-G_t C_t+(B_t-G_tD_t)H_t \\ \end{array}\right]\\
\tilde B_t&= \left[\begin{array}{cc} B_t  \\    B_t-G_tD_t \\ \end{array}\right]  L_t ,\;
\tilde C_t = \left[\begin{array}{cc} C_t &  D_t H_t \\ \end{array}\right],\;
\tilde  D_t = D_t L_t \nn} and the input is WGN with unit variance. Note that, to construct the least favorable model,
first we generate the gains $G_t$ performing a forward sweep of the robust filter (\ref{rob_est_formula1})-(\ref{rob_est_formula2}) over interval $[0,T]$, then we generate the  matrices $\Omega_{t}$
through a backward sweep over interval $[0,T]$. Therefore, increasing the simulation interval beyond $[0,T]$ requires performing a new backward sweep of recursion (\ref{recursion_in_omega}). 
Then, we can evaluate the performance of an arbitrary estimator \alg{\label{arb_filter}\hat x^\prime_{t+1}=A_t \hat x_t^\prime +G_t^\prime (y_t-C_t \hat x_t^\prime)}
 applied to the the least favorable model. Let
\alg{\Pi_t=\Es\left[\left[\begin{array}{c} e_t^\prime \\ e_t  \\ \end{array}\right] \left[\begin{array}{cc} (e_t^\prime)^T & e_t^T \\\end{array}\right]\right] \nn}
where $e_t$ is the estimation error of the optimal filter (\ref{rob_est_formula1}) and $e_t^\prime$ is the estimation error of filter (\ref{arb_filter}). Then, it can be proven
that $\Pi_t$ obeys to the {\em Lyapunov} equation, \cite{ROBUST_STATE_SPACE_LEVY_NIKOUKHAH_2013},
\alg{\label{iter_lyap}\Pi_{t+1}=& \left( \tilde A_t- \left[\begin{array}{c} G_t^\prime \\0  \\ \end{array}\right] \tilde C_t\right)\Pi_t \left( \tilde A_t- \left[\begin{array}{c} G_t^\prime \\0  \\ \end{array}\right] \tilde C_t\right)^T\nn\\ &+ \left( \tilde B_t- \left[\begin{array}{c} G_t^\prime \\0  \\ \end{array}\right] \tilde D_t\right)\left( \tilde B_t- \left[\begin{array}{c} G_t^\prime \\0  \\ \end{array}\right] \tilde D_t\right)^T }
where $\Pi_0=I_2\otimes V_0$.

\section{Risk Sensitive Filtering} \label{sec_RS}

Consider the robust Kalman filter we presented in Section \ref{sec_rob_kalman} with $\tau=0$. If we replace $\theta_t$ with a constant value $\theta$ we recognize immediately that
we obtain the risk sensitive filter, \cite{OPTIMAL_SPEYER_FAN_BANAVAR_1992,RISK_PROP_BANAVAR_SPEIER_1998,RISK_WHITTLE_1980}. This suggest us
that the risk sensitive filter can extended using the $\tau$-divergence family.
Consider the {\em Markov-Gauss} state space model (\ref{state_space_model}). Let $\bar{f}_t(z_t|Y_{t-1})$ be the conditional density of $z_t$
given $Y_{t-1}$ based on the model  (\ref{state_space_model}) and defined in (\ref{bar_f_t_Gaussian}).
The classic risk sensitive estimator $g_t^{\circ}$
at time $t$ is defined as
\alg{\label{RS_est} g^{\circ}_t= \underset{g_t\in\Gc_t}{\mathrm{argmin}}\Es_{\bar{f}_t}[\exp(\theta \|x_{t+1}-g_t(y_t)\|^2)\,|\,Y_{t-1}]}
where $\Gc_t$ is the set of estimators for which the objective function in (\ref{RS_est}) is finite. $\theta>0$ is the risk sensitivity parameter. More precisely, the larger $\theta$ is the more the objective function in (\ref{RS_est})  penalizes estimators with
large errors.  In \cite{boel2002robustness}, it has been shown that the risk sensitive estimator is solution to the following minimax problem 
\alg{\small \label{RS_est_minimax}g_t^{\circ}= \underset{g_t\in\Gc_t}{\mathrm{argmin}}  \underset{\tilde f_t\in\Bc_t}{\max}\Es_{\tilde f_t}[ & \|x_{t+1}  -g_t(y_t)\|^2 \,|\,Y_{t-1}]\nn\\ & -\theta^{-1}\Ds_{KL}(\tilde f_t\| f_t)}
where $\Bc_t=\{\tilde f_t \hbox{ s.t. } \Ds_{KL}(\tilde f_t \| f_t)<\infty  \}$.
The second term in the objective function in (\ref{RS_est_minimax}) is always nonpositive because $\Ds_{KL}(\tilde f_t\| f_t)\geq 0$.
For small values of $\theta$, it takes large negative values for conditional densities not close to the nominal one. Therefore, the maximizer is obliged to choose a conditional density close to the nominal one.
On the contrary, for large values of $\theta$,   it takes (negative) values close to zero for some conditional densities not close to the nominal one. In such  a situation, the maximizer has the possibility to  choose those conditional densities.
Note that, this behaviour  does not change if we replace $\Ds_{KL}$ with another divergence measure.

In our setting  $\bar{f}_t(z_t|Y_{t-1})$ is Gaussian by assumption. In addition, if we assume that $\tilde f_t(z_t|Y_{t-1})$ is Gaussian, then  $\Ds_{KL}(\tilde f_t \| f_t)=\Dc_0(\tilde f_t \| f_t)$
where $\Dc_0$ has been defined in (\ref{def_tau_div}). It is then natural to extend the minimax problem (\ref{RS_est_minimax}) to the $\tau$-divergence family:
 \alg{g_t^{\circ}=\underset{g_t\in\Gc_t}{\mathrm{argmin}}  \underset{\tilde f_t\in\Bc_t}{\max}\Es_{\tilde f_t}[ & \|x_{t+1}  -g_t(y_t)\|^2 \,|\,Y_{t-1}]\nn\\ & -\theta^{-1}\Dc_{\tau}(\tilde f_t\| f_t) \nn}

By applying Corollary \ref{corollario_static}, the optimal $\tau$-risk sensitive estimator takes the form of (\ref{rob_est_formula1})-(\ref{rob_est_formula2})
where
  \alg{ \label{Vt1_rs}& V_{t+1}=\nn\\
&\left\{
      \begin{array}{ll}
        L_{P_{t+1}}\left(I_n-\theta(1-\tau) L_{P_{t+1}}^TL_{P_{t+1}}\right)^{\frac{1}{\tau-1}} L_{P_{t+1}}^T , & \hspace{-0.1cm}  0<\tau<1\\
L_{P_{t+1}}\exp\left(\theta L_{P_{t+1}}^T L_{P_{t+1}}\right) L_{P_{t+1}}^T, & \tau=1.
      \end{array}
    \right.} It is worth noting that, for the case $0 \leq \tau< 1$, $V_{t+1}$ is defined provided that $0<P_{t+1}<(\theta(1-\tau))^{-1}I_n$, while for the case $\tau=1$,
    it is well defined whenever $P_{t+1}$ is positive definite. The algorithmic scheme is similar to Algorithm \ref{algorithm}: the unique difference is that Step 4 is now removed. Finally, while the risk sensitivity parameter of the  robust filter of Section \ref{sec_rob_kalman} is time varying, and its evolution is governed by $c_t$, now it is constant.

\section{Simulation Results} \label{sec_simu}
\begin{figure}[htbp]
\begin{center}
\includegraphics[width=\columnwidth]{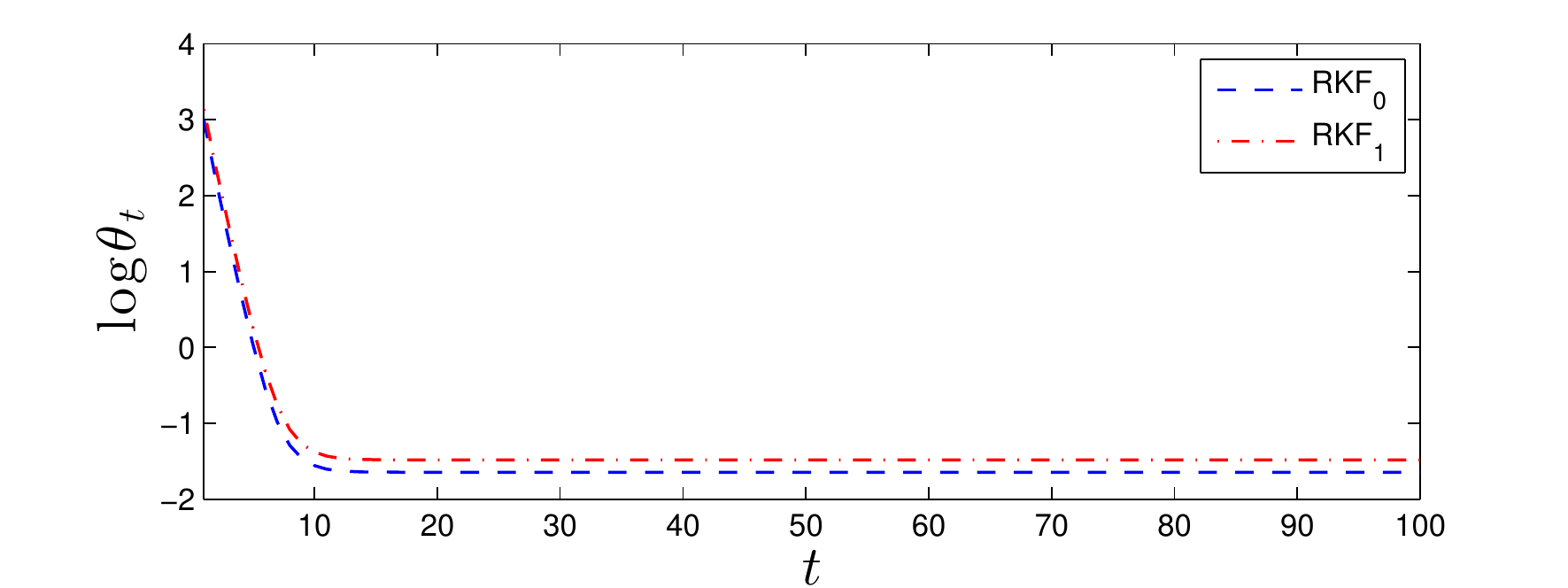}
\end{center}
 \caption{Evolution of the risk sensitive parameter $\theta_t$ with $c=10^{-1}$.}\label{fig3}
\end{figure}

We consider the time-invariant model (\ref{state_space_model}) with 
\alg{ A&=\left[\begin{array}{cc}0.1 & 1  \\0 & 1.2 \\\end{array}\right], & B=\left[\begin{array}{ccc} 0.01 & 0 & 0  \\0 & 0.01 &0 \\\end{array}\right]\nn\\
C&=\left[\begin{array}{cc} 1 & -1  \\ \end{array}\right], &D=\left[\begin{array}{ccc} 0 & 0  & 0.1  \\ \end{array}\right]\nn }
and $x_0\sim \Nc(0,V_0)$ with $V_0=0.01\cdot I_2$. 
We consider the following three filters:
KF is the standard Kalman filter; RKF$_0$ is the robust Kalman filter of Section \ref{sec_rob_kalman} with  $\tau=0$;  RKF$_1$ is the robust Kalman filter of Section \ref{sec_rob_kalman} with  $\tau=1$. In Figure \ref{fig3} we show the evolution of the risk sensitivity parameter of RKF$_0$ and RKF$_1$ for $c=10^{-1}$. We notice the $\theta_t$ is constant in the steady state, that is RKF$_0$ and RKF$_1$ coincide with the risk sensitive filters of Section \ref{sec_RS} with $\tau=0$, $\theta\approx 0.19$ and $\tau=1$, $\theta\approx 0.23$, respectively, in the steady state.    

In what follows we evaluate the performance of RK, RKF$_0$ and RKF$_1$, which have structure (\ref{arb_filter}), applied to the least favorable model (\ref{ls_model}). More precisely, for each filter, applied to (\ref{ls_model}), we consider the estimation error $e^\prime_t=[\,e_t^1\, e^2_t\,]^T$. Then, we compute the variance of $e_t^1$ and $e_t^2$ through (\ref{iter_lyap}). We consider two situations: $c$ large, i.e. nominal and least favorable model are very different; $c$ small, i.e. nominal and least favorable model are similar.

\subsection{Large tolerance}
\begin{figure}[htbp]
\begin{center}
\includegraphics[width=\columnwidth]{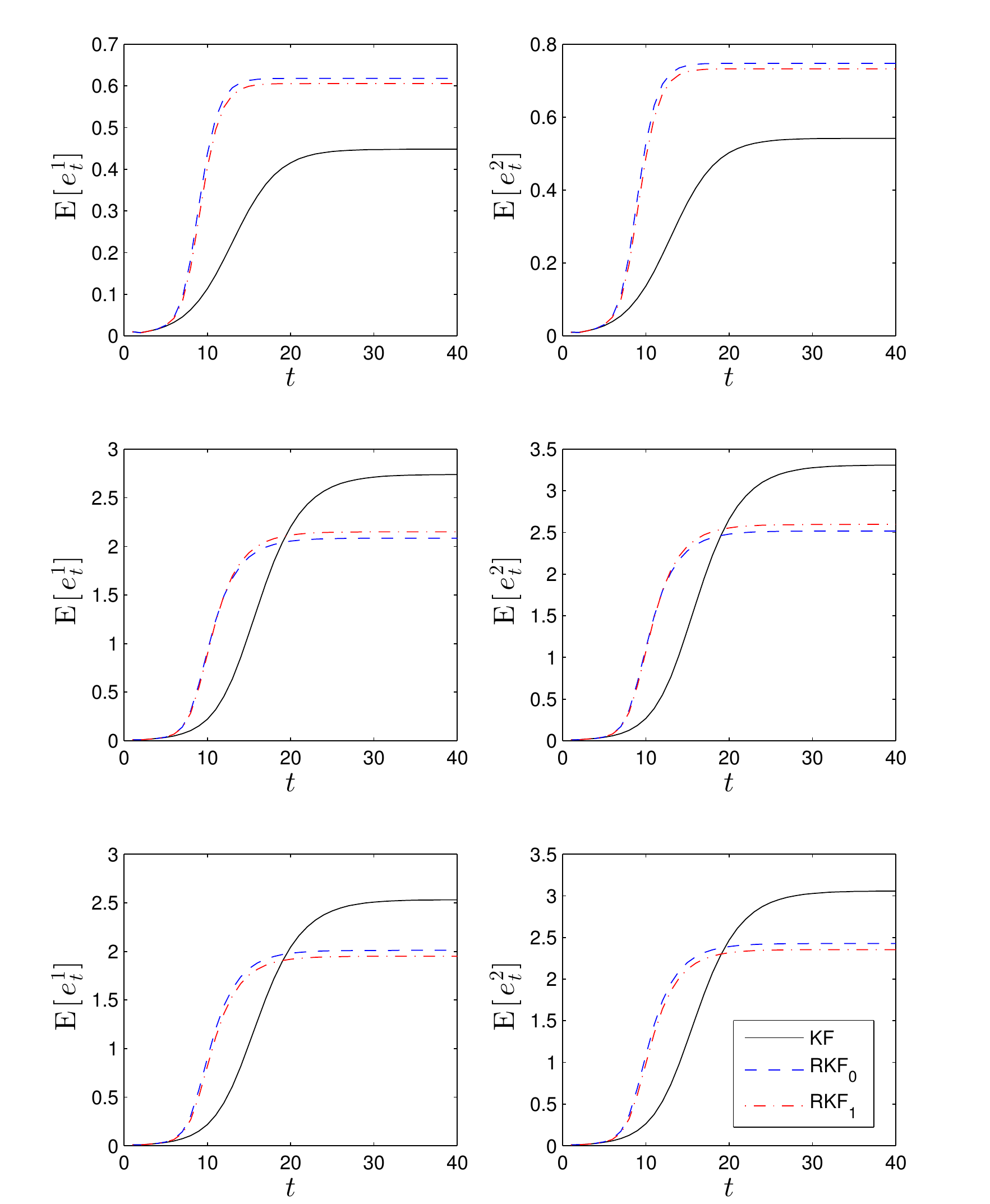}
\end{center}
 \caption{Variances of $e_t^1$ and $e_t^2$ when the filters are applied to the nominal model (first row); to (\ref{ls_model})
 with $\tau=0$ and $c=10^{-1}$ (second row); to (\ref{ls_model}) with $\tau=1$ and $c=10^{-1}$ (third row). Here RKF$_0$ and RKF$_1$ have $c=10^{-1}$.}\label{fig1}
\end{figure}

Here RKF$_0$ and RKF$_1$ have tolerance $c=10^{-1}$.
In the first experiment, we apply these filters to the nominal model (\ref{state_space_model}). The variance of $e_t^1$ and $e_t^2$  are depicted in the first row of Figure \ref{fig1}. As expected, KF performs better than the others. Moreover, the variances of RKF$_0$ are slightly larger than the ones of RKF$_1$.
In the second experiment, we apply these filters to the 
least favorable model (\ref{ls_model}) with $\tau=0$ and $c=10^{-1}$. The variances of  $e_t^1$ and $e_t^2$  are depicted in the second row of Figure \ref{fig1}. Obviously, RKF$_0$ is the best estimator because it has been designed 
with respect to this model. Although RKF$_1$ has been designed with respect to another model, it performs better than KF.
In the third experiment, we apply these filters to the 
least favorable model (\ref{ls_model}) with $\tau=1$ and $c=10^{-1}$.  The variances of $e_t^1$ and $e_t^2$ are depicted in the third row of Figure \ref{fig1}. In this case RKF$_1$ is the best estimator because it is optimal with respect to the underlying model. Also in this case, the worst estimator is KF. 
From these experiments we can conclude that:
\begin{itemize}
\item the smaller $\tau$ is, the more conservative the filter is, that is, the smaller $\tau$ is, the more the uncertainty class contains models with larger mean square error. This property has been noticed  also for the static estimation problem in \cite{OPTIMALITY_ZORZI}.
\item the family of robust Kalman filters provide better performances than the standard Kalman filter, even in the case that the least favorable model belongs to an uncertainty class parametrized by a different $\tau$. 
 \end{itemize}

\subsection{Small tolerance}
\begin{figure}[htbp]
\begin{center}
\includegraphics[width=\columnwidth]{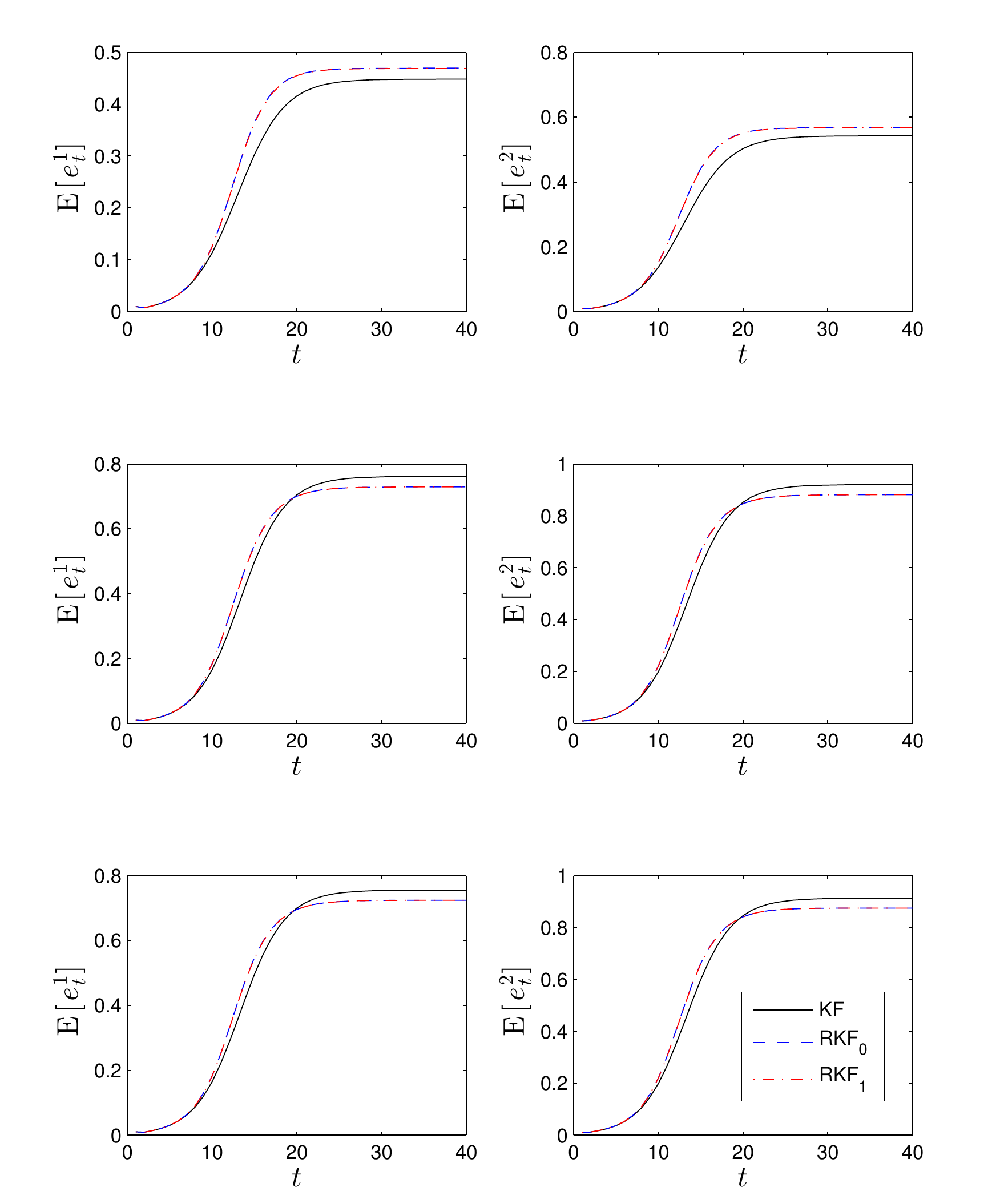}
\end{center}
 \caption{Variances of $e_t^1$ and $e_t^2$ when the filters are applied to the nominal model (first row); to (\ref{ls_model})
 with $\tau=0$ and $c=5\cdot 10^{-3}$ (second row); to (\ref{ls_model}) with $\tau=1$ and $c=5\cdot 10^{-3}$ (third row). Here RKF$_0$ and RKF$_1$ have $c=5\cdot 10^{-3}$.}\label{fig2}
\end{figure}

We perform the same experiment three experiments of before where the unique difference is the tolerance which now is $c=5\cdot 10^{-3}$, see Figure \ref{fig2}. RKF$_0$ and RKF$_1$ provides the same performance which is comparable with the one of KF. Therefore, as long as the discrepancy between the nominal and the least favorable model is not too large, then the performance of KF applied to (\ref{ls_model})
 does not deteriorate too much. 

\section{Conclusions} \label{sec_concl}
In this paper, we have considered a robust filtering problem under incremental model perturbations
characterized by the $\tau$-divergence family. 
The family of robust estimators we proposed is the solution to a minimax problem. These robust estimators have an iterative structure similar to the one of the Kalman filter. We have derived the corresponding least favorable models. 
Moreover, we have extended the risk sensitive filter to a family of risk sensitive like filters. 
Finally, a simulation study shows that parameter $\tau$ tunes how conservative the robust filter is.

\section*{acknowledgments}

I am very grateful and indebted to Prof. Bernard C. Levy. During my visit at UC Davis in 2011, he has been always willing to clear up my doubts arising 
from his papers  \cite{LEVY_NIKOUKHAH_2004,ROBUST_STATE_SPACE_LEVY_NIKOUKHAH_2013}.


\end{document}